\documentclass{amsart}

\usepackage{amssymb}
\usepackage{picture}

\newtheorem{theorem}{Theorem}[section]
\newtheorem{lemma}[theorem]{Lemma}
\newtheorem{example}[theorem]{Example}
\newtheorem{proposition}[theorem]{Proposition}

\newtheorem{corollary}[theorem]{Corollary}

\newcommand{\N}{\mathbb N}

\newcommand{\R}{\mathbb R}

\newcommand{\Z}{\mathbb Z}

\newcommand{\on}{\operatorname}
\author{Artur Bartoszewicz}
\address{Institute of Mathematics, \L \'od\'z University of Technology,
W\'olcza\'nska 215, 93-005 \L \'od\'z, Poland}
\email {arturbar@p.lodz.pl}

\author{Szymon G\l \c ab}
\address{Institute of Mathematics, \L \'od\'z University of Technology,
W\'olcza\'nska 215, 93-005 \L \'od\'z, Poland}
\email {szymon.glab@p.lodz.pl}

\author{Jacek Marchwicki}
\address{Institute of Mathematics, \L \'od\'z University of Technology,
W\'olcza\'nska 215, 93-005 \L \'od\'z, Poland}
\email {marchewajaclaw@gmail.com}

\title[Achievement sets of conditionally convergent series]{Achievement sets of conditionally convergent series}

\thanks{The second author has been supported by the National Science Centre Poland Grant no. DEC-2012/07/D/ST1/02087.}

\begin{document}

\begin{abstract}
Considering the sets of subsums of series (or achievement sets) we show that for conditionally convergent series the multidimensional case is much more complicated than that of the real line. Although we are far from the full topological classification of such sets, we present many surprising examples and catch the ideas standing behind them in general theorems. 
 
\end{abstract}

\maketitle

\section{Introduction}

Probably, S. Kakeya \cite{Kakeya} was the first one to consider topological properties of subsums of absolutely convergent series of real numbers. For the absolutely summable sequence $(x_{n})$, we call the set $\on{A}(x_{n})=\{\sum_{n=1}^\infty\varepsilon_nx_n : (\varepsilon_{n})\in\{0,1\}^{\mathbb{N}}\}$ \emph{the set of subsums} or \emph{the achievement set} \cite{Jones}. Of course, for the sequence $(x_{n})$ with almost all terms equal to zero, the set $\on{A}(x_{n})$ is finite. Kakeya has shown:
\begin{theorem}
 For an absolutely summable sequence with infinitely many nonzero terms 
\begin{itemize}
\item The set $\on{A}(x_{n})$ is a compact perfect set.
\item If for almost all $n$ $\vert x_{n}\vert >\sum_{k>n}\vert x_{k}\vert$ then $\on{A}(x_{n})$ is homeomorphic to the Cantor set (after M. Moran we call such sequences \emph{quickly convergent})
\item  If for almost all $n$ $\vert x_{n}\vert \leq\sum_{k>n}\vert x_{k}\vert$ then $\on{A}(x_{n})$ is a finite union of closed intervals. Moreover it can be reversed for nonincreasing sequences  $(\vert x_{n}\vert)$.
\end{itemize}
\end{theorem}
Kakeya conjectured that Cantor-like sets and finite unions of closed intervals are the only possible achievement sets for sequences $(x_{n})\in \ell_{1}\setminus c_{00}$. The Kakeya results were rediscovered many times and his conjecture was repeated, even when the first counterexamples were given. The first counterexamples were published by Weinstein and Shapiro \cite{WS}, Ferens \cite{Ferens} and Guthrie and Nymann \cite{GN}. Due to Guthrie, Nymann and Saenz \cite{GN,NS1} we know that the achievement set of an absolutely summable sequence can be a finite set, a finite union of intervals, homeomorphic to the Cantor set or it can be a so called Cantorval. A Cantorval is a set homeomorphic to the union of the Cantor set and sets which are removed from the unite segment by even steps of the Cantor set construction. The last result gives the partition of $\ell_{1}$ into four disjoint sets. Topological and algebraic properties of these sets were recently considered in \cite{BBGS,BBGS1}. Some sufficient conditions for a given sequence to be a Cantorval were recently described in \cite{BBFS,BFS,Jones}. The connections between achievement sets of some absolutely summable sequences and self-similar sets were observed in the last of the above papers. 

If $\sum_{n=1}^\infty x_n$ is an absolutely convergent series in Banach space, then the function $\{0,1\}^\N\ni(\varepsilon_n)\mapsto\sum_{n=1}^\infty\varepsilon_nx_n$ is a continuous function which maps the Cantor space $\{0,1\}^\N$ onto $\on{A}(x_n)$, see for example \cite{BG}. In particular, the achievement set $\on{A}(x_n)$ is compact. We will prove that this function is continuous also if the series $\sum_{n=1}^\infty x_n$ is unconditionally convergent. It will follow immediately from Lemma \ref{UnconditionLemma}. It is well-known \cite[Theorem 1.3.2, p. 10]{KadKad} that a series $\sum_{n=1}^\infty x_n$ in a Banach space is unconditionally convergent if and only if each series of the form $\sum_{n\in A}x_n$, $A\subseteq\N$, is convergent. 
\begin{lemma}\label{UnconditionLemma}
Assume that $\sum_{n=1}^\infty x_n$ is an unconditionally convergent series in a Banach space $X$. Then for every $\varepsilon>0$ there is $N\in\N$ such that for every set $A\subseteq\N$
$$
\Vert\sum_{k\in A\setminus\{1,\dots,N\}}x_k\Vert\leq\varepsilon.
$$
\end{lemma}
\begin{proof}
Suppose to the contrary that there is $\varepsilon>0$ such that for every $N$ there is $A\subseteq\{N+1,N+2,\dots\}$ with
$$
\Vert\sum_{k\in A}x_k\Vert>\varepsilon.
$$
For $N=1$ find $A_1$ with $\Vert\sum_{k\in A_1}x_k\Vert>\varepsilon$. There is a finite set $F_1\subseteq A_1$ with $\Vert\sum_{k\in F_1}x_k\Vert>\varepsilon$. In the second step for $N=\max F_1$ find $A_2\subseteq\{N+1,N+2,\dots\}$ with $\Vert\sum_{k\in A_2}x_k\Vert>\varepsilon$. As before we take a finite set $F_2\subseteq A_2$ with $\Vert\sum_{k\in F_2}x_k\Vert>\varepsilon$. Proceeding inductively, we produce finite sets $F_1,F_2,\dots$ such that $\max F_i<\min F_{i+1}$ and $\Vert\sum_{k\in F_i}x_k\Vert>\varepsilon$. Put $A=\bigcup_{i\geq 1}F_i$. Then by the Cauchy condition the series $\sum_{n\in A}x_n$ diverges which yields a contradiction.
\end{proof}

One can define the achievement sets also for all sequences in Banach spaces. Then there should be considered only these sequences $(\varepsilon_{n})\in\{0,1\}^{\mathbb{N}}$ for which a series $\sum_{n=1}^\infty\varepsilon_n x_n$ is convergent. 
\begin{theorem}
For sequences of reals with $\lim\limits_{n\to\infty} x_n=0$ we have:
\begin{itemize}
\item A series $\sum_{n=1}^\infty x_n$ is potentially conditionally convergent (both series of positive and negative terms are divergent) if and only if $\on{A}(x_{n})=\mathbb{R}$. 
\item A series of negative terms is convergent and a series of positive terms is divergent (or vice versa) if and only if the achievement set of $(x_n)$ is a half line.
\end{itemize}
\end{theorem}
For simple proofs see for example \cite{BFPW,Jones,N1}. If $(x_n)$ does not converge to zero, then $\on{A}(x_{n})$ is always $F_{\sigma}$-set \cite{Jones}. So, for conditionally convergent series of reals the achievement set $\on{A}(x_{n})$ is equal to $\mathbb{R}$, which is exactly the same set as the sum range $\on{SR}(x_{n})$, i.e. the set of rearrangements $\sum_{n=1}^\infty x_{\sigma(n)}$,  by the classical Riemann Theorem on permutations of conditionally convergent series. The aim of the present paper is to show that a situation in multidimensional spaces is quite different. We observe among others that for the achievement set  $\on{A}(x_{n})$ of conditionally convergent series in $\mathbb{R}^2$ the following are possible. 
\begin{itemize}
\item The intersection of $\on{A}(x_{n})$  and $\on{SR}(x_{n})$ could be a singleton and moreover we mention that it is always nonempty set (Example \ref{InstructiveExample});
\item $\on{A}(x_{n})$ can be a graph of function (Example \ref{ExampleAchievementIsGraph}); 
\item $\on{A}(x_{n})$ can be a dense set in $\mathbb{R}^2$ with an empty interior (Example \ref{ExampleOfSecondTypeSeries}); 
\item $\on{A}(x_{n})$ can be neither $F_{\sigma}$ nor $G_{\delta}$-set (Theorem \ref{NeitherFSigmaNorGDelta});
\item $\on{A}(x_{n})$ can be an open set not equal to the whole $\mathbb{R}^2$ (Theorem \ref{OpenAchievementSet}); 
\end{itemize}
On the other hand to obtain regular achievement sets let us make a simple observation. Let $X$ be a Banach space. Assume that $A=\on{A}(x_n)$ and $B=\on{A}(y_n)$ are achievement sets in $X$. Then
$$
A\times B=\on{A}((x_1,0),\:(0,y_1),\:(x_2,0),\:(0,y_2),...).
$$ 
If $T:X\to Y$ is a bounded linear operator from $X$ to some other Banach space $Y$, then
$$
T(A)=\on{A}(Tx_1,Tx_2,Tx_3,\dots).
$$
Now take any conditionally convergent series $\sum_{n=1}^\infty x_n$ and absolutely convergent series $\sum_{n=1}^\infty y_n$, both on the real line. Then $\on{A}(x_{n})=\R$ and $\on{A}(y_{n})=C$ is a compact set. By the above observation there are conditionally convergent series on the plane such that their achievement sets equal $\R\times C$, $\R^2$ or any rotation of $\R\times C$. 

By our knowledge the achievement sets of series in multidimensional spaces were considered in a few papers only. For example in \cite{M1,M2} the author has studied quickly convergent series. In \cite{BG} the series considered by the authors are also absolutely convergent. So our paper is probably the first one on achievement sets of conditionally convergent series in $\mathbb{R}^{n}$ for $n$ greater than one. On the other hand properties of sum range sets are well-studied \cite{KadKad}. The strong suggestion to consider achievement sets of series in multidimensional spaces was given by Nitecki in his nice lecture \cite{N1} (the shorter version of this survey is \cite{N2}). 

\section{Cardinality of achievement sets}

In this section we study the cardinality of achievement sets in Banach spaces. Then we will prove that the achievement set of conditionally convergent set is \emph{perfectly dense-in-itself}, that is for any point $x$ which can be achieved and any $\varepsilon>0$, the intersection of the achievement set and the ball $B(x,\varepsilon)$ contains a perfect set. 

\begin{proposition}\label{PropCardinality}
Let $X$ Banach space and let $(x_n)$ be a sequence of elements of $X$. Then  
\begin{itemize}
\item[(i)] $\on{A}(x_{n})$ is finite, if there are finite many non-zero $x_{n}$'s;
\item[(ii)] $\on{A}(x_{n})$ is infinite and countable, if there are infinite many non-zero $x_{n}$'s and there is $\delta>0$, such that $\Vert x_{k}\Vert\geq\delta$ for any nonzero element $x_k$; moreover if $X$ is finitely dimensional, then $\on{A}(x_{n})$ is unbounded;
\item[(iii)] $\on{A}(x_{n})$ contains a perfect set otherwise.
\end{itemize}
\end{proposition}
\begin{proof}
(i) Let $x_{k_{1}},\ldots,x_{k_{m}}$ be a list of all non-zero terms of $(x_{n})$. Then the cardinality of $\on{A}(x_{n})$ is not greater than $2^{m}$.

(ii) Since zero terms do not affect the achievement set, we may assume that our sequence $(x_{n})$ consists of non-zero elements. Firstly, we consider the case of $X=\mathbb{R}^{k}$ with the supremum norm. Then $x_n=(x_n(1),\dots,x_n(k))$. Without lose of generality we may assume that the set $F_{j}=\{n\in\mathbb{N} : x_{n}(j)\geq\delta\}$ is infinite for some $j\leq k$. Hence $\sum_{n\in F_{j}} x_{n}(j)=\infty$, so the achievement set is unbounded and therefore infinite. 

In general, if the set $\{x_{n}:n\in\N\}$ spans a finitely dimensional space, then we may assume that $X$ is finitely dimensional and therefore isomorphic to $\mathbb{R}^{k}$ with the supremum norm. If $\{x_{n}:n\in\N\}$ spans an infinitely dimensional space, then $\{x_{n}:n\in\N\}$ is an infinite set contained in $\on{A}(x_{n})$.  

To see that $\on{A}(x_{n})$ is countable, note that no subsequence of $x_{n}$ converges to zero. Therefore every element of $\on{A}(x_{n})$ is a sum of finitely many $x_n$'s. 

(iii) The negation of the first two conditions means that there exists an infinite subsequence $(x_{n_{l}})$ of non-zero terms, which tends to $0$. We may assume $\Vert x_{n_{l+1}}\Vert<\frac{\Vert x_{n_{l}}\Vert}{3}$ for every $l\in\mathbb{N}$. Then $\{0,1\}^\N\ni(\varepsilon_{l})\stackrel{f}{\mapsto}\sum_{l=1}^{\infty}\varepsilon_{l}x_{n_{l}}$ is injective. To see this assume that $\overline{\varepsilon_{j}}\neq\varepsilon_{j}$ for some $j\in\N$. We have 
$$
\Vert\sum\limits_{l=j+1}^{\infty}\varepsilon_{l}x_{n_{l}}-\sum\limits_{l=j+1}^{\infty}\overline{\varepsilon}_{l}x_{n_{l}}\Vert=\Vert\sum\limits_{l=j+1}^{\infty}(\varepsilon_{l}-\overline{\varepsilon}_{l})x_{n_{l}}\Vert\leq\sum\limits_{l=j+1}^{\infty}\Vert x_{n_{l}}\Vert\leq\frac{3}{2}\Vert x_{n_{j+1}}\Vert<\Vert x_{n_{j}}\Vert.
$$
Hence 
$$
\Vert\sum\limits_{l=j}^{\infty}\varepsilon_{l}x_{n_{l}}-\sum\limits_{l=j}^{\infty}\overline{\varepsilon}_{l}x_{n_{l}}\Vert=
\Vert\varepsilon_{j}x_{n_{j}}-\overline{\varepsilon_{j}}x_{n_{j}}+\sum\limits_{l=j+1}^{\infty}\varepsilon_{l}x_{n_{l}}-\sum\limits_{l=j+1}^{\infty}\overline{\varepsilon}_{l}x_{n_{l}}\Vert\geq
$$
$$
\geq\Vert\varepsilon_{j}x_{n_{j}}-\overline{\varepsilon_{j}}x_{n_{j}}\Vert-\Vert\sum\limits_{l=j+1}^{\infty}\varepsilon_{l}x_{n_{l}}-\sum\limits_{l=j+1}^{\infty}\overline{\varepsilon}_{l}x_{n_{l}}\Vert=
\Vert x_{n_{j}}\Vert-\Vert\sum\limits_{l=j+1}^{\infty}\varepsilon_{l}x_{n_{l}}-\sum\limits_{l=j+1}^{\infty}\overline{\varepsilon}_{l}x_{n_{l}}\Vert>0.
$$ 
This implies that there exists $r<j$ such that $\overline{\varepsilon_{r}}\neq\varepsilon_{r}$. After finitely many steps we get  $\overline{\varepsilon_{1}}\neq\varepsilon_{1}$. But $$\Vert\sum\limits_{l=1}^{\infty}\varepsilon_{l}x_{n_{l}}-\sum\limits_{l=1}^{\infty}\overline{\varepsilon}_{l}x_{n_{l}}\Vert\geq\Vert x_{n_{1}}\Vert-\Vert\sum\limits_{l=2}^{\infty}\varepsilon_{l}x_{n_{l}}-\sum\limits_{l=2}^{\infty}\overline{\varepsilon}_{l}x_{n_{l}}\Vert>0.
$$
Since $\sum_{l=1}^\infty x_{n_l}$ is absolutely convergent, the mapping $f$ is continuous. Hence $\on{A}(x_{n_{l}})\subseteq\on{A}(x_n)$ is a continuous injective image of the Cantor space $\{0,1\}^\N$. 
\end{proof}

\begin{proposition}
Let $\sum_{n=1}^\infty x_n$ be a conditionally convergent series in a Banach space $X$. Then $\on{A}(x_n)$ is perfectly dense in itself. 
\end{proposition}

\begin{proof}
Note that the set $\{\sum_{n\in F}x_n:F$ is finite$\}$ is dense in $\on{A}(x_n)$. Let $x\in\on{A}(x_n)$ and $\varepsilon>0$. There is a finite set $F$ such that $\Vert\sum_{n\in F}x_n-x\Vert<\varepsilon/2$. Using the same method as in the proof of Proposition \ref{PropCardinality}(iii), we find a subsequence $x_{n_l}$ of $x_n$ such that $\{n_l:l\in\N\}\cap F=\emptyset$, $\sum_{l=1}^\infty\Vert x_{n_l}\Vert<\varepsilon/2$ and $\on{A}(x_{n_l})$ is a perfect set. Then $\sum_{n\in F}x_n+\on{A}(x_{n_l})\subseteq B(x,\varepsilon)\cap\on{A}(x_n)$. 
\end{proof}

\section{Achievement sets of conditionally convergent series}

We start this section from giving an instructive example of conditionally convergent series $\sum_{n=1}^\infty(\frac{(-1)^{n+1}}{n},\frac{1}{2^n})$ on the plane. The achievement set of this series have several properties which show that the theory of achievement sets of conditionally convergent series in multidimensional spaces is much more complicated and interesting than that of one-dimensional series. The analysis of this example will lead us to make some general observations:
\begin{itemize}
\item First we note in Proposition \ref{RozkladCondConv} that every conditionally convergent series in $\R^n$ either has a sum range equal to $\R^n$ (\emph{first type series}) or it is, up to linear isometry, of the form $\sum_{n=1}^\infty(x_n,y_n)$ where $x_n\in\R^k$, $y_n\in\R^{m-k}$, $\sum_{n=1}^\infty x_n$ is conditionally convergent with a sum range equal to $\R^k$ and $\sum_{n=1}^\infty y_n$ is absolutely convergent (\emph{second type series}).
\item We will observe that the closure of $\on{A}(\frac{(-1)^{n+1}}{n},\frac{1}{2^n})$ contains the sum range $\on{SR}(\frac{(-1)^{n+1}}{n},\frac{1}{2^n})$. This is true in every Banach space (Lemma \ref{AIsClosedInSR}).
\item We will observe that the closure of $\on{A}(\frac{(-1)^{n+1}}{n},\frac{1}{2^n})$ equals $\R\times\on{A}(\frac{1}{2^n})$. A similar fact is true in Euclidean spaces (Theorem \ref{AIsDenseInProduct}).
\item The achievement set $\on{A}(\frac{(-1)^{n+1}}{n},\frac{1}{2^n})$ is not closed. This phenomena is generalized in Theorem \ref{TheoremAIsNotClosed}.
\item A series $\sum_{n=1}^\infty(\frac{2}{3^n},\frac{(-1)^{n+1}}{n})$ considered in Example \ref{ExampleAchievementIsGraph}, a slight modification of that from Example \ref{InstructiveExample}, has the achievement set which is neither $F_\sigma$ nor $G_\delta$. A wide class of series with that property is given in Theorem \ref{NeitherFSigmaNorGDelta}. In the proof of Theorem \ref{NeitherFSigmaNorGDelta} we observe that an achievement set of conditionally convergent series is always an analytic (or $\Sigma^1_1$) set.    
\end{itemize}
The series from Example \ref{InstructiveExample} is of the first type. By Lemma \ref{AIsClosedInSR} the achievement set of series of the second type is dense in the whole space. One can very easily give an example of series which achievement set is actually the whole space. It is much harder to give an example of series of the second type which achievement is smaller. We will present it in Example \ref{ExampleOfSecondTypeSeries} -- its achievement set will be a null subset of the plane. 

The classical theorem of Steinitz states that if $\sum_{n=1}^\infty x_n$ is conditionally convergent in $\R^m$, then the set $\on{SR}(x_n)=\{\sum_{n=1}^\infty x_{\sigma(n)}:\sigma\in S_\infty\}$ of all convergent rearrangements $\sum_{n=1}^\infty x_{\sigma(n)}$ of $\sum_{n=1}^\infty x_n$ is an affine subspace of $\R^m$. More precisely, if $\Gamma=\{f\in(\R^m)^*:\sum_{n=1}^\infty \vert f(x_n)\vert<\infty\}$ is a set of all convergent functionals on $\sum_{n=1}^\infty x_{n}$, and $\Gamma^\perp=\{x\in\R^m:f(x)=0$ for all $f\in\Gamma\}$ is an anihilator of $\Gamma$, then 
$$
\on{SR}(x_n)=\sum_{n=1}^\infty x_n+\Gamma^\perp.
$$ 

\begin{example}\label{InstructiveExample}
\emph{
Let $x_{n}=(x_{n}^{(1)},x_{n}^{(2)})=(\frac{(-1)^{n+1}}{n},\frac{1}{2^n})$ be terms of a series $\sum _{n=1}^{\infty} x_{n}$ in $\mathbb{R}^{2}$. Then clearly $\sum _{n=1}^{\infty} x_{n}=(\log 2,1)$. By Riemann theorem for every $x\in\mathbb{R}$ one can find $\sigma\in S_{\infty}$ such that $x=\sum _{n=1}^{\infty} x_{\sigma(n)}^{(1)}$. Since the permutation of indices does not affect the value of absolutely convergent series, we have $\sum _{n=1}^{\infty} x_{\sigma(n)}^{(2)}=1$. Hence $\on{SR}(x_{n})=\mathbb{R}\times\{1\}$. Let $D=\{\sum _{n=1}^{k}\frac{\varepsilon_{n}}{2^{n}}: (\varepsilon_{n})\in\{0,1\}^{k}, k\in\mathbb{N}\}$ be the set of all dyadic numbers from interval $[0,1)$. Then for every $d\in D$ there are $k\in\mathbb{N}$ and $(\varepsilon_{n})\in\{0,1\}^{k}$ with $d=\sum _{n=1}^{k}\frac{\varepsilon_{n}}{2^{n}}$. Put $F_{d}:=\{n\leq k : \varepsilon_{n}=1\}$. After removing finitely many terms from conditionally convergent series, we still have conditionally convergent series.  Therefore $\on{SR}((x_n)_{n\in\mathbb{N}\setminus{F_{d}}})=\mathbb{R}\times\{1-d\}$. 
}

\emph{
From Lemma \ref{AIsClosedInSR} we get $\on{SR}((x_n)_{n\in\mathbb{N}\setminus{F_{d}}})\subseteq\overline{\on{A}((x_n)_{n\in\mathbb{N}\setminus{F_{d}}})}\subseteq \overline{\on{A}(x_{n})}$. Since $d\in D\iff 1-d\in D$, then $\bigcup _{d\in D} (\mathbb{R}\times\{d\})\subseteq \overline{\on{A}(x_{n})}$ and consequently $\overline{\bigcup _{d\in D} (\mathbb{R}\times\{d\}})\subseteq \overline{\on{A}(x_{n})}$. But $D$ is dense in the interval $[0,1]$, so we have:
$$
\overline{\bigcup _{d\in D} (\mathbb{R}\times\{d\}})=\mathbb{R}\times\overline{\bigcup _{d\in D}\{d\}}=\mathbb{R}\times\overline{D}=\mathbb{R}\times [0,1].
$$
Thus $\mathbb{R}\times [0,1]\subseteq \overline{\on{A}(x_{n})}$. The reverse inclusion is obvious; therefore $\overline{\on{A}(x_{n})}=\mathbb{R}\times [0,1]$. Suppose $(z,1)\in \on{A}(x_{n})$. The only way to get $\sum _{n=1}^{\infty}\varepsilon_{n} x_{n}^{(2)}=1$ for some $(\varepsilon_{n})\subseteq \{0,1\}^{\mathbb{N}}$ is to take $\varepsilon_{n}=1$ for each $n\in\mathbb{N}$. Hence $z=\log 2$. That proves that the only point which belongs to an achievement set with the second coordinate $1$ is $(\log 2,1)$. Thus $\on{A}(x_{n})$ is not closed. Finally note that $\on{SR}(x_n)=\R\times\{1\}$. Therefore $\on{A}(x_n)\cap \on{SR}(x_n)=\{(\log 2,1)\}$. 
}
\end{example}

The following observation, which is probably mathematical folklore, allows us to consider only two types of conditionally convergent series in Euclidean spaces. 

\begin{proposition}\label{RozkladCondConv}
Let $\sum_{n=1}^{\infty} x_{n}$ be a conditionally convergent series in  $\mathbb{R}^m$ such that $\dim(\Gamma^\perp)=k<m$. Then there exists an isomorphism $(T_1,T_2)\in L(\mathbb{R}^m,\mathbb{R}^k\times \mathbb{R}^{m-k})$ such that $\sum_{n=1}^{\infty} T_{1}(x_{n})$ is conditionally convergent with $SR(\sum_{n=1}^{\infty} T_{1}(x_{n}))=\mathbb{R}^k$ and $\sum_{n=1}^{\infty} T_{2}(x_{n})$ is absolutely convergent in $\mathbb{R}^{m-k}$.
\end{proposition}

\begin{proof} 
Let $k$ be a dimension of $\Gamma^\perp$ and let $Y$ be orthogonal to $\Gamma^\perp$ so that $\R^m=\Gamma^\perp\oplus Y$. Let $e_1,\dots,e_m$ be a standard basis of $\R^m$ and let $e_1',\dots,e_m'$ be an orthogonal basis of $\Gamma^\perp\oplus Y$ such that $\Gamma^\perp=\on{span}\{e_1'\dots,e_k'\}$ and $Y=\on{span}\{e_{k+1}'\dots,e_{m}\}$. Let $T:\R^m\to\Gamma^\perp\oplus Y$ be a linear isomorphism such that $T(e_i)=e_i'$ for every $i=1,\dots,m$. For $x=\sum_{i=1}^mx(i)e_i$ let $T_1(x)=\sum_{i=1}^kx(i)e_i'$ and $T_2(x)=\sum_{i=k+1}^mx(i)e_i'$. Then $T=(T_1,T_2)$. Let $f_i(x)=x(i)$. Then $f_i\in\Gamma$ for $i>k$, which means that $\sum_{n=1}^\infty\Vert T_2(x_n)\Vert<\infty$. 

Let $\Lambda=\{f\in(\Gamma^\perp)^*:\sum_{n=1}^\infty \vert f(T_1(x_n))\vert<\infty\}$. Let $\pi_{\leq k}:\R^m\to\R^k$ be a projection onto the first $k$ coordinates. For $f\in(\Gamma^\perp)^*$ define $\tilde{f}\in(\R^m)^*$ as $\tilde{f}(x)=f(\pi_{\leq k}(x))$. Then $f\in\Lambda\iff \tilde{f}\in\Gamma\iff\tilde{f}=0\iff f=0$. Thus $\Lambda=\{0\}$ and by Stenitz theorem $\on{SR}(T_1(x_n))=\Lambda^\perp=\R^k$. 
\end{proof}

Note that $T$ from Proposition \ref{RozkladCondConv} does not change geometrical and topological properties of subsets of $\R^m$. Therefore we will assume that a conditionally convergent series in $\R^m$ can be rearranged to get any point of $\R^m$ or it is of the form $\sum_{n=1}^\infty(x_n,y_n)$ where $x_n\in\R^k$, $y_n\in\R^{m-k}$, $\on{SR}(x_n)=\R^k$ and $\sum_{n=1}^\infty y_n$ is absolutely convergent.

The following lemma shows a relation between the sum range and the achievement set of the series.
\begin{lemma}\label{AIsClosedInSR}
Let $\sum _{n=1}^{\infty} x_{n}$ be a conditionally convergent series in a Banach space $X$. Then $\on{SR}(x_{n})\subseteq\overline{\on{A}(x_{n})}$.
\end{lemma} 

\begin{proof}
 Let $\varepsilon>0$ and $x\in\on{SR}(x_{n})$, then $x=\sum _{n=1}^{\infty} x_{\sigma(n)}$ for some $\sigma\in S_{\infty}$. One can find natural number $k$ such that $\Vert x-\sum _{n=1}^{k} x_{\sigma(n)}\Vert<\varepsilon$. Denote $A=\{m : \sigma(n)=m, n\leq k\}=\sigma(\{1,\ldots,k\})$ and define the sequence $\varepsilon_{n}=1$ for $n\in A$ and $\varepsilon_{n}=0$ otherwise. Then  $\sum _{n=1}^{k} x_{\sigma(n)}=\sum _{n=1}^{\infty} \varepsilon_{n} x_{n}$, so $\Vert x-\sum _{n=1}^{\infty} \varepsilon_{n} x_{n}\Vert<\varepsilon$. Hence $x\in\overline{\on{A}(x_{n})}$.
\end{proof}

Lemma \ref{AIsClosedInSR} implies that the achievement sets of conditionally convergent series in finite dimensional spaces are unbounded. The situation changes in infinitely dimensional Banach spaces. 
\begin{example}\label{C0Example} 
\emph{
Let $(e_n)$ be a standard basis of $c_0$. Let $(x_n)$ be a sequence of the form
$$
e_1,-e_1,\frac12 e_2,-\frac12 e_2,\frac12 e_2,-\frac12 e_2,\frac13 e_3,-\frac13 e_3,\frac13 e_3,-\frac13 e_3,\frac13 e_3,-\frac13 e_3,\dots
$$
Then $\sum_{n=1}^\infty x_n$ is convergent to zero, while its rearrangement
$$
e_1-e_1+\frac12 e_2+\frac12 e_2-\frac12 e_2-\frac12 e_2+\frac13 e_3+\frac13 e_3+\frac13 e_3-\frac13 e_3-\frac13 e_3-\frac13 e_3+\dots
$$ 
is divergent, since the sequence of its partial sums contains each $e_i$. Since the projection of the series on each coordinate contains only finitely many nonzero terms and a finite sum does not change under rearrangements, $\on{SR}(\sum_{n=1}^\infty x_n)=\{0\}$. Let $X_n=\{\frac{k}{n}:k\in[-n,n]\cap\Z\}$. Note that $\on{A}(\sum_{n=1}^\infty x_n)=(\prod_{n=1}^\infty X_n)\cap c_0$. Therefore the achievement set $\on{A}(\sum_{n=1}^\infty x_n)$ is closed and bounded. 
}
\end{example}

It can happen that $\on{SR}(x_{n})=\on{A}(x_{n})$. For example \begin{itemize}
\item $\on{A}(\frac{(-1)^n}{n},0)=\R\times\{0\}=\on{SR}(\frac{(-1)^n}{n},0)$;
\item $\on{A}(\frac{(-1)^n}{n},\frac{(-1)^n}{n})=\{(x,x)\in\R\}=\on{SR}(\frac{(-1)^n}{n},\frac{(-1)^n}{n})$;
\item Let the series $\sum_{n=1}^\infty x_n$ and $\sum_{n=1}^\infty y_n$ be conditionally convergent on the real line. Then 
$$\on{A}((x_1,0),(0,y_1),(x_2,0),(0,y_2),\dots)=\R^2=\on{SR}((x_1,0),(0,y_1),(x_2,0),(0,y_2),\dots).$$
\end{itemize}

Note that if $\sum _{n=1}^{\infty} x_{n}$ is a conditionally convergent series in a Banach space $X$, then the intersection of $\on{A}(x_n)$ and $\on{SR}(x_n)$ is non empty; more precisely $\sum _{n=1}^{\infty} x_{n}\in \on{A}(x_n)\cap \on{SR}(x_n)$. Note that in Example \ref{InstructiveExample} the intersection $\on{A}(x_n)\cap \on{SR}(x_n)$ is actually a singleton. 

Next example is a slight modification of Example \ref{InstructiveExample}. 

\begin{example}\label{ExampleAchievementIsGraph}
\emph{
Let $(x_{k},y_k)=(\frac{2}{3^k},\frac{(-1)^{k+1}}{k})$. Since the mapping $\{0,1\}^\N\ni(\varepsilon_{k})\mapsto\sum_{k=1}^{\infty}\varepsilon_{k}x_{k}$ is one-to-one, the achievement set $\on{A}(x_{k},y_{k})$ is a graph of function $f$ with a domain contained in the ternary Cantor set. Moreover $f$ maps onto $\mathbb{R}$, and its domain is $F_{\sigma\delta}$ set which is not $G_{\delta\sigma}$, which will follow from Theorem \ref{NeitherFSigmaNorGDelta}.
}
\end{example}

Using the idea from Example \ref{InstructiveExample} we will prove the following. 
\begin{theorem}\label{AIsDenseInProduct}
Let $(x_{n})$ be a sequence in $\mathbb{R}^k$ such that $\sum _{n=1}^{\infty} x_{n}$ is conditionally convergent with $\on{SR}(x_{n})= \mathbb{R}^k$ and let $(y_{n})$ be a sequence in $\mathbb{R}^m$ such that a series $\sum _{n=1}^{\infty} y_{n}$ is absolutely convergent. Then 
$\overline{\on{A}(x_{n},y_{n})}=\mathbb{R}^k\times \on{A}(y_{n})$.
\end{theorem}
\begin{proof} ''$\subseteq$''. It is easy to see that $\on{A}(x_{n},y_{n})\subseteq\on{A}(x_{n})\times\on{A}(y_{n})$. Therefore
$$
\overline{\on{A}(x_{n},y_{n})}\subseteq\overline{\on{A}(x_{n})\times  \on{A}(y_{n})}= \overline{\on{A}(x_{n})}\times\overline{ \on{A}(y_{n})}.
$$
Since $\on{SR}(x_{n})= \mathbb{R}^k$, by Lemma \ref{AIsClosedInSR} we get  $\overline{\on{A}(x_{n})}=\R^k$. By the absolute convergence of $\sum _{n=1}^{\infty} y_{n}$ we get the compactness of its achievement set, so  $\overline{\on{A}(y_{n})}= \on{A}(y_{n})$. 

''$\supseteq$''. Let $(x,y)\in\mathbb{R}^k\times \on{A}(y_{n})$ and $\varepsilon>0$. Since $y=\sum _{n=1}^{\infty} \varepsilon_{n} y_{n}$ for some  $(\varepsilon_{n})\in\{0,1\}^{\mathbb{N}}$, there exists $k_{\varepsilon}\in\mathbb{N}$ such that $\Vert y-\sum _{n=1}^{N}\varepsilon_{n} y_{n}\Vert<\varepsilon$ for every $N\geq k_{\varepsilon}$.  Since $\sum _{n=1}^{\infty} y_{n}$ is absolutely convergent we may assume that $\sum _{n=k_{\varepsilon}+1}^{\infty}\Vert y_{n}\Vert<\varepsilon$. Let $K=\{n\leq k_{\varepsilon} : \varepsilon_{n}=1\}=\{k_{1}<\ldots<k_{l}\}$. Define $\sigma(n)=k_{n}$ for $n\in\{1,\ldots,l\}$. 

Note that if $\on{SR}(x_{n})$ is the whole space $\mathbb{R}^k$, then $\on{SR}((x_n)_{n\geq k_{\varepsilon}+1})$ is the whole space $\mathbb{R}^k$ as well. In particular 
$x-\sum _{n=1}^{l}x_{k_n}\in\on{SR}((x_n)_{n\geq k_{\varepsilon}+1})$. Therefore
there exists $M\in\mathbb{N}$ and a one-to-one mapping $\tau:\{k_{\varepsilon}+1,\ldots,M\}\rightarrow\{k_{\varepsilon}+1,k_{\varepsilon}+2,\ldots\}$ such that $\Vert x-\sum _{n=1}^{l}x_{k_n}-\sum _{n=k_{\varepsilon}+1}^{M}x_{\tau(n)}\Vert<\varepsilon$. Enumerate the range $\tau(\{k_{\varepsilon}+1,k_{\varepsilon}+2,\ldots,M\})$ as $\{k_{l+1}<\ldots<k_{l+l'}\}$ and define the sequence: 
\begin{displaymath}
\delta_{n} = \left\{ \begin{array}{ll} \varepsilon_{n}, &  n\leq k_{\varepsilon}; \\ 1, &  n=k_{l+i} \ \ \text{for} \ \ i\in\{1,\ldots,l'\}; \\ 0, & \textrm{otherwise} .\end{array} \right.
\end{displaymath} 
Therefore
$$\Vert \sum _{n=1}^{\infty} \varepsilon_{n} y_{n}-\sum _{n=1}^{\infty} \delta_{n} y_{n}\Vert=\Vert \sum _{n=k_{\varepsilon}+1}^{\infty} \varepsilon_{n} y_{n}-\sum _{n=k_{\varepsilon}+1}^{\infty} \delta_{n} y_{n}\Vert=\Vert \sum _{n=k_{\varepsilon}+1}^{\infty} (\varepsilon_{n}-\delta_{n}) y_{n}\Vert\leq\sum _{n=k_{\varepsilon}+1}^{\infty}\Vert y_{n}\Vert<\varepsilon. $$
Hence $\Vert y-\sum _{n=1}^{\infty} \delta_{n} y_{n}\Vert<\varepsilon$. Moreover 
$$\Vert x-\sum _{n=1}^{\infty} \delta_{n} x_{n}\Vert=\Vert x-\sum _{n\leq k_{\varepsilon}} \delta_{n} x_{n}-\sum _{n=1}^{l'} \delta_{k_{l+n}} x_{k_{l+n}}\Vert=\Vert x-\sum _{n\leq k_{\varepsilon}} \varepsilon_{n} x_{n}-\sum _{n=l+1}^{l+l'}  x_{k_{n}}\Vert=$$ 
$$\Vert x-\sum _{n\in K}  x_{n}-\sum _{n=l+1}^{l+l'}  x_{k_{n}}\Vert=\Vert x-\sum _{n=1}^{l}  x_{k_{n}}-\sum _{n=k_{\varepsilon}+1}^{M}  x_{\tau(n)}\Vert<\varepsilon.$$
Since the norms in finite dimensional spaces are equivalent, there is $C>0$ such that 
$$\Vert (x,y)-\sum _{n=1}^{\infty}\delta_{n}(x_{n},y_{n})\Vert\leq C\cdot\max\{\Vert x-\sum _{n=1}^{\infty} \delta_{n} x_{n}\Vert,\Vert y-\sum _{n=1}^{\infty} \delta_{n} y_{n}\Vert\}<C\cdot\varepsilon.$$
Finally we get $(x,y)\in \overline{\on{A}(x_{n},y_{n})}$ which finishes the proof.
\end{proof}

Now, we present the sufficient condition for conditionally convergent series to have a not closed achievement set. Here we use again the idea from Example \ref{InstructiveExample}. 

\begin{lemma}\label{UniqueRepresentationOfExtremePoint}
Let $\sum_{n=1}^{\infty} y_{n}$ be an absolutely convergent series in $ \mathbb{R}^m$ with $y_{n}\neq 0$ for each $n\in\mathbb{N}$. Then for every extreme point $a$ of the achievement set $\on{A}(y_{n})$, there is a unique sequence $(\varepsilon_{n})\in\{0,1\}^{\mathbb{N}}$ such that $a=\sum_{n=1}^{\infty} \varepsilon_{n}y_{n}$.
\end{lemma}
\begin{proof} We will show that $a$ is achieved for a unique sequence $(\varepsilon_{n})\in\{0,1\}^{\mathbb{N}}$. Suppose in contrary that $a=\sum_{n=1}^{\infty} \varepsilon_{n}y_{n}=\sum_{n=1}^{\infty} \delta_{n}y_{n}$ for two distinct sequences $(\varepsilon_n)$ and $(\delta_n)$,  and put $M=\{n\in\mathbb{N}: \varepsilon_{n}\neq\delta_{n}\}$. Divide $M$ into two disjoint sets $M_{\varepsilon}=\{n\in M : \varepsilon_{n}=1, \delta_{n}=0\}$ and  $M_{\delta}=\{n\in M : \varepsilon_{n}=0, \delta_{n}=1\}$. Then $a=\sum_{n\in M_{\varepsilon}}  y_{n} + \sum_{n\in M^{c}}\varepsilon_{n} y_{n} = \sum_{n\in M_{\delta}}  y_{n} + \sum_{n\in M^{c}} \varepsilon_{n} y_{n} $, so 
$$a=\frac{1}{2} \sum\limits_{n\in M^{c}}\varepsilon_{n} y_{n}+\frac{1}{2} (\sum_{n\in M} y_{n}+\sum_{n\in M^{c}} \varepsilon_{n} y_{n})=\frac{1}{2}b+\frac{1}{2}c.$$
Since $b= \sum_{n\in M^{c}}\varepsilon_{n} y_{n}$ and $c=\sum_{n\in M} y_{n}+\sum_{n\in M^{c}} \varepsilon_{n} y_{n}$ we have $b,c\in\on{A}(y_{n})$. Since $a$ is the extreme point of $\on{A}(y_{n})$, we get $b=c$, so $\sum_{n\in M} y_{n}=0$. Since each $y_n$ is nonzero, $M$ has at least two elements. Assume that $n_{0}\in M$, then $y_{n_{0}}+\sum_{l\in M\setminus\{n_{0}\}} y_{l}=0$. Define $b'=b+y_{n_{0}}$ and $c'=c-y_{n_{0}}$. Then $b',c'\in\on{A}(y_{n})$. We have $a=\frac{1}{2} b'+\frac{1}{2} c'$, and using again the assumption that $a$ is the extreme point, we get $b'=c'$. But then $y_{n_{0}}=\sum_{l\in M\setminus\{n_{0}\}} y_{l}$ and $\sum_{n\in M} y_{n}=0$, and therefore $y_{n_{0}}=0$, which gives us a contradiction. 
\end{proof}
\begin{theorem}\label{TheoremAIsNotClosed}
Let  $\sum_{n=1}^{\infty} x_{n}$ be a conditionally convergent series in $\mathbb{R}^k$ with $\on{SR}(x_{n})= \mathbb{R}^k$ and $\sum_{n=1}^{\infty} y_{n}$ be an absolutely convergent series with $y_{n}\neq 0$ for each $n\in\mathbb{N}$. Then $\on{A}(x_{n},y_{n})$ is not  closed.
\end{theorem}
\begin{proof}
Let $a$ be an extreme point of the achievement set $\on{A}(y_{n})$; such a point exists because the achievement set $\on{A}(y_{n})$ is compact. By Lemma \ref{UniqueRepresentationOfExtremePoint} a point $a$ has a unique representation $a=\sum_{n=1}^{\infty}\varepsilon_{n}y_{n}$. Hence, $a$-section of $\on{A}(x_{n},y_{n})$ is a singleton if $\sum_{n=1}^\infty\varepsilon_n x_n$ converges, or it is empty if $\sum_{n=1}^\infty\varepsilon_n x_n$ diverges. On the other hand $\{(x,a) : x\in\R^{k}\}\subseteq\R^{k}\times \on{A}(y_{n})=\overline{\on{A}(x_{n},y_{n})}$. Hence, $\on{A}(x_{n},y_{n})$ is not closed.
\end{proof}

Now, we will prove that the achievement set of conditionally convergent series does not need to be neither $F_\sigma$ nor $G_\delta$. This will apply to the series from Example \ref{ExampleAchievementIsGraph}. 

Let us start from giving definitions of some important notions from the descriptive set theory. A topological space is Polish if it is completely metrizable and separable. An $F_{\sigma\delta}$ subset $A$ of a Polish space $X$ is called $\Pi^0_3$-complete if for any zero-dimensional Polish space $Y$ and for any $F_{\sigma\delta}$ subset $B$ of $Y$ there is a continuous function $f:Y\to X$ such that $f^{-1}(A)=B$. It is known that $\Pi^0_3$-complete sets are $F_{\sigma\delta}$ but not $G_{\delta\sigma}$ sets. To prove that an $F_{\sigma\delta}$ subset $C$ of Polish space $Z$ is $\Pi^0_3$-complete, it is enough to take a known example $A$ of a $\Pi^0_3$-complete subset of a Polish space $X$ and find a continuous function $g:X\to Z$ such that $g^{-1}(C)=A$. For more information we refer the reader to \cite{Kechris}. 

\begin{proposition}\label{Pi03CompleteSet}
Let $\sum_{n=1}^\infty x_n$ be a conditionally convergent series in $\R$ such that $\on{A}(x_n)=\R$. 
Then the set $E:=\{(\varepsilon_n)\in 2^\N:\sum_{n=1}^\infty\varepsilon_nx_n$ converges$\}$ is $\Pi^0_3$-complete subset of $\{0,1\}^\N$. 
\end{proposition} 

\begin{proof}
Note that 
$$
(\varepsilon_n)\in E\iff\forall k\in\N\;\exists l\in\N\;\forall M>m\geq l\;\;\;\left(\vert\sum_{n=m}^{M}\varepsilon_nx_n\vert\leq\frac{1}{k}\right).
$$
Therefore $E$ is $F_{\sigma\delta}$ subset of $2^\N$. 

To prove that $E$ is $\Pi^0_3$-complete we will use the fact that the set 
$$
C_3:=\{v\in\N^\N:\lim_{n\to\infty}v(n)=\infty\}
$$
is $\Pi^0_3$-complete, for details see \cite[Section 23A]{Kechris}. It is enough to construct a continuous function $\psi:\N^\N\to\{0,1\}^\N$ such that $v\in C_3\iff\psi(v)\in E$. Spaces $\{0,1\}^\N$ of 0-1 sequences and $\N^\N$ of sequences of natural numbers are considered with the metric $d(x,y)=2^{-n}$ where $n=\min\{k:x(k)\neq y(k)\}$.

One can define inductively sets $F_n=F_n(v)$ and $H_n=H_n(v)$ such that $F_0=H_0=\emptyset$ and for every $n\geq 1$
\begin{itemize}
\item[(i)] $F_n<H_n<F_{n+1}$, that is $\max F_n<\min H_n$, $\max H_n<\min F_n$;
\item[(ii)] $\vert\sum\limits_{k\in F_1\cup H_1\cup\dots\cup F_n}x_k\vert<2^{-n}$;
\item[(iii)] $\vert\sum\limits_{k\in F_1\cup H_1\cup\dots\cup F_n\cup H_n}x_k-2^{-v(n)}\vert<2^{-n}$;
\item[(iv)] $x_k<0$ for $k\in\bigcup_{n\geq 1}F_n$ and $x_k>0$ for $k\in\bigcup_{n\geq 1}H_n$.
\end{itemize}
Note that the above construction can be made uniformly in the sense that if $v(i)=v'(i)$ for $i\leq n$, then $F_i(v)=F_i(v')$ and $H_i(v)=H_i(v')$ for $i\leq n$.   
 
Now, we define $\psi:\N^\N\to\{0,1\}^\N$ as follows. Let $\psi(v)$ be a characteristic function of $\bigcup_{n=1}^\infty F_n(v)\cup H_n(v)$. Since the construction is uniform, then for $v,v'\in\N^\N$ such that $v(i)=v'(i)$ for $i\leq n$ we have $d(\psi(v),\psi(v'))\leq 2^{-n}$. Therefore $\psi$ is continuous. We will prove that $v\in C_3\iff\sum_{n=1}^\infty\psi(v)(n)x_n$ is convergent. 

If $v\notin C_3$, then there are $m\in\N$ and an infinite set $L\subseteq\N$ such that $v(l)=m$ for all $l\in L$. Thus by the construction, a series $\sum_{n=1}^\infty\psi(v)(n)x_n$ diverges, since the sequence of partial sums  $(\sum_{n=1}^N\psi(v)(n)x_n)$  has, by (ii) and (iii), two accumulation points - $0$ and $2^{-m}$.

If $v\in C_3$, then $2^{-v(n)}\to 0$. Let $\varepsilon>0$. There is $l\in\N$ such that $2^{-v(n)},2^{-n}<\varepsilon/2$ for $n\geq l$. Let $n\geq l$. By (ii) and (iii) we obtain
$$
\vert\sum\limits_{k\in F_1\cup H_1\cup\dots\cup F_n}x_k\vert=\vert\sum_{k=1}^{\max F_n}\psi(v)(k)x_k\vert<2^{-n}<\varepsilon
$$ 
and 
$$
\vert\sum\limits_{k\in F_1\cup H_1\cup\dots\cup F_n\cup H_n}x_k-2^{-v(n)}\vert=
\vert\sum_{k=1}^{\max H_n}\psi(v)(k)x_k-2^{-v(n)}\vert<2^{-n}+2^{-v(n)}<\varepsilon.
$$
By (iv) we have also
$$
\vert\sum_{k=1}^{M}\psi(v)(k)x_k\vert<2^{-n}
$$
for every $M\geq\max H_l$. Thus $\sum_{n=1}^\infty\psi(v)(n)x_n$ converges to zero. 
\end{proof}

The similar result to Proposition \ref{Pi03CompleteSet} was proved by Cohen in \cite{C} -- having conditionally convergent series $\sum_{n=1}^\infty x_n$ in $\R^n$, the set of all permutations $\sigma$ such that $\sum_{n=1}^\infty x_{\sigma(n)}$ converges, is a $\Pi^0_3$-complete subset of $S_\infty$.  The constructed example of a $\Pi^0_3$-complete set in Proposition \ref{Pi03CompleteSet} will be used to prove the following.

\begin{theorem}\label{NeitherFSigmaNorGDelta}
Let $\sum_{n=1}^\infty x_n$ be a conditionally convergent series in $\R$, and let $\sum_{n=1}^\infty y_n$ be an absolutely convergent series such that the function $\{0,1\}^\N\ni(\varepsilon_n)\mapsto\sum_{n=1}^\infty\varepsilon_n y_n$ is one-to-one. Then the achievement set $\on{A}(x_n,y_n)$ is a Borel subset of $\R^2$ which is neither $G_\delta$ nor $F_\sigma$. 
\end{theorem}

\begin{proof}
Note that 
$$
(x,y)\in\on{A}(x_n,y_n)\iff\exists(\varepsilon_n)\in 2^\N\;\forall m\;\exists l\;\forall k\geq l\; (\vert x-\sum_{n=1}^k\varepsilon_nx_n\vert<\frac1m\text{ and }\vert y-\sum_{n=1}^k\varepsilon_ny_n\vert<\frac1m)
$$
which shows that $\on{A}(x_n,y_n)$ is $\Sigma^1_1$ as a projection of Borel set, and 
$$
(x,y)\in\on{A}(x_n,y_n)\iff\exists!(\varepsilon_n)\in 2^\N\;\forall m\;\exists l\;\forall k\geq l\; (\vert x-\sum_{n=1}^k\varepsilon_nx_n\vert<\frac1m\text{ and }\vert y-\sum_{n=1}^k\varepsilon_ny_n\vert<\frac1m)
$$
which shows that $\on{A}(x_n,y_n)$ is $\Pi^1_1$, see \cite[18.11]{Kechris}. By Suslin Theorem $\on{A}(x_n,y_n)$ is Borel. 

By Theorem \ref{AIsDenseInProduct} the achievement set $\on{A}(x_n,y_n)$ is dense in $\R\times\on{A}(y_n)$. However every horizontal section of $\on{A}(x_n,y_n)$ consists of at most one point. Suppose that $\on{A}(x_n,y_n)$ is a $G_\delta$. Then it would be a comeager subset of $\R\times\on{A}(y_n)$, and therefore almost all of its horizontal sections in the sense of category would be comeager in $\R$. That would be a contradiction since each its horizontal section is at most one-point. 

Suppose that $\on{A}(x_n,y_n)$ is an $F_\sigma$ subset of $\R^2$. Then it would be a countable union of compact sets, and consequently its projection $\on{proj}_2(\on{A}(x_n,y_n))$ on the second coordinate would be a countable union of compact sets, that is an $F_\sigma$ set. But $\on{proj}_2(\on{A}(x_n,y_n))$ is homeomorphic to the set $E$ from Lemma \ref{Pi03CompleteSet} which is not $F_\sigma$. A contradiction. 
\end{proof}

We have also considered whether or not in Theorem \ref{AIsDenseInProduct} the condition that the sum range of the conditionally convergent series is the whole space, can be replaced by the same condition for the achievement set of the series. The following proposition shows us that in some cases we can reverse  Theorem \ref{AIsDenseInProduct}.
\begin{proposition}\label{PropAIsDenseSoIsSR}
Assume that $\sum_{n=1}^{\infty} x_{n}$ is conditionally convergent in $\mathbb{R}^{k}$. If $\on{A}(x_{n})$ is a dense subset of $\mathbb{R}^{k}$, then  $\on{SR}(x_{n})=\mathbb{R}^{k}$.
\end{proposition}
\begin{proof}
Suppose that $\on{SR}(x_{n})\neq\mathbb{R}^{k}$. By Proposition \ref{RozkladCondConv} we have $\sum_{n=1}^{\infty} x_{n}=\sum_{n=1}^{\infty}T(y_{n},z_{n})$, where $T: \mathbb{R}^{m}\times  \mathbb{R}^{k-m}\rightarrow  \mathbb{R}^{k}$ is an isomorphism and $\on{SR}(y_{n})=\mathbb{R}^{m}$ for some $m$ with $1\leq m<k$ and $\sum_{n=1}^{\infty} z_{n}$ is an absolutely convergent series in $\mathbb{R}^{k-m}$. Without loss of generality we may assume $x_{n}=(y_{n},z_{n})$. Hence $\sum_{n=1}^{\infty} x_{n}=\sum_{n=1}^{\infty}(y_{n},z_{n})$. Thus by Theorem \ref{AIsDenseInProduct} and compactness of $\on{A}(z_{n})$ we have $\overline{\on{A}(x_{n})}=\mathbb{R}^{m}\times \on{A}(z_{n})\neq \mathbb{R}^{k}$ which yields a contradiction.
\end{proof}

By Theorem \ref{AIsDenseInProduct} and Proposition \ref{PropAIsDenseSoIsSR} we immediately obtain the following.

\begin{corollary}\label{CorAIsDenseInProduct}
Let $(x_{n})\subseteq \mathbb{R}^k$ be such that $\sum _{n=1}^{\infty} x_{n}$ is conditionally convergent and such that $\on{A}(x_{n})$ is dense in $\mathbb{R}^k$. Let $(y_{n})\subseteq \mathbb{R}^m$ be such that a series $\sum _{n=1}^{\infty} y_{n}$ is absolutely convergent. Then 
$\overline{\on{A}(x_{n},y_{n})}=\mathbb{R}^k\times \on{A}(y_{n})$.
\end{corollary}

On the other hand, one can construct the series on the plane, which sum range is $\mathbb{R}^{2}$ but the achievement set is a dense set of measure zero. 
\begin{example}\label{ExampleOfSecondTypeSeries}
\emph{
Let  $\sum_{n=1}^{\infty} (x_{n},y_{n})$ be defined as $x_{i}=\frac{(-1)^{i}}{2^{10^{k^{2}}}}$ and $y_{i}=\frac{(-1)^{i}}{2^{k}}$ for $i\in (n_{k-1},n_{k}]$, where $n_{0}=0$ and $n_{k+1}=n_{k}+2^{10^{k^{2}}+1}$.
Since $\sum_{n=1}^{\infty} (x_{n},y_{n})$ is alternating and $ (x_{n},y_{n})\rightarrow (0,0)$, the series $\sum_{n=1}^{\infty} (x_{n},y_{n})$ is convergent to $(0,0)$. Firstly, we will show that $\on{SR}(x_{n},y_{n})=\mathbb{R}^{2}$. Let $\Gamma=\{f\in (\mathbb{R}^{2})^{*}: \sum_{n=1}^{\infty} \vert f(x_{n},y_{n})\vert<\infty\}$ be a set of all convergence functionals for our series. We need to prove that $\Gamma$ does not contain nontrivial functionals. Every $f\in (\mathbb{R}^{2})^{*}$ is of the form $f(x,y)=ax+by$ for some $a,b\in\mathbb{R}$. We have 
$$ \sum_{i=n_{k-1}+1}^{n_{k}} \vert y_{i}\vert=(n_{k}-n_{k-1})\frac{1}{2^{k}}=2^{10^{k^{2}}+1-k}\text{ and } \sum_{i=n_{k-1}+1}^{n_{k}} \vert x_{i}\vert=(n_{k}-n_{k-1})\frac{1}{2^{10^{k^{2}}}}=2.$$
If $a=0$ and $b\neq 0$, then 
$$\sum_{i=1}^{\infty} \vert f(x_{i},y_{i})\vert=\sum_{k=1}^{\infty}\sum_{i=n_{k-1}+1}^{n_{k}} \vert b y_{i}\vert=\vert b\vert\sum_{k=1}^{\infty}2^{10^{k^{2}}+1-k}=\infty.$$
If  $a\neq 0$ and $b=0$, then 
$$\sum_{i=1}^{\infty} \vert f(x_{i},y_{i})\vert=\sum_{k=1}^{\infty}\sum_{i=n_{k-1}+1}^{n_{k}} \vert a x_{i}\vert=\vert a\vert\sum_{k=1}^{\infty}2=\infty.$$
Now let us consider $a\neq 0$ and $b\neq 0$. See that for large a enough natural number $k$, which satisfies the inequality $\frac{\vert a\vert}{2^{10^{k^{2}}-k}}<\frac{\vert b\vert}{2}$, we have $\vert\frac{a}{2^{10^{k^{2}}}}+\frac{b}{2^{k}}\vert\geq\frac{\vert b\vert}{2^{k}}-\frac{\vert a\vert}{2^{10^{k^{2}}}}\geq\frac{\vert b\vert}{2^{k+1}}$ and the series $\sum_{k=1}^{\infty}\sum_{i=n_{k-1}+1}^{n_{k}} \frac{\vert b\vert}{2^{k+1}}$  diverges. Consequently
$$\sum_{i=1}^{\infty} \vert f(x_{i},y_{i})\vert=\sum_{k=1}^{\infty}\sum_{i=n_{k-1}+1}^{n_{k}} \vert a x_{i}+b y_{i}\vert=\infty.$$
From the Steiniz theorem we have $\on{SR}(x_{n},y_{n})=\mathbb{R}^{2}$.
}

\emph{
Now we will show that  $\on{A}(x_{n},y_{n})\subseteq(L\cup \mathbb{Q})\times\mathbb{R}$, where $L=\{x : \forall{r\in\mathbb{N}}\;\exists{p,q\in\mathbb{Z}}\;(0<\vert x-\frac{p}{q}\vert<\frac{1}{q^{r}})\}$ is the set of all Liouville numbers on the real line. It is well-known that $L$ has the Lebesgue measure zero. 
Let $r\in\mathbb{N}$. Suppose that $(x,y)\in \on{A}(x_{n},y_{n})$, that is $(x,y)=\sum_{n=1}^{\infty}\varepsilon_{n}(x_{n},y_{n})$. Then there exists $l\in\mathbb{N}$ such that $\vert\sum_{n=N}^{M} \varepsilon_{n} y_{n}\vert\leq 1$ for every $M>N\geq l$.
There exists $k_{0}\in\mathbb{N}$ for which $n_{k_{0}-1}$ is not less than $l$, so $\vert\sum_{i=n_{k-1}+1}^{n_{k}} \varepsilon_{i} y_{i}\vert\leq 1$ for every $k\geq k_{0}$. Note that $\lim_{k\to\infty}\frac{10^{k^{2}}-k-1}{10^{(k-1)^{2}}}=\infty$ and the sequence $(\frac{10^{k^{2}}-k-1}{10^{(k-1)^{2}}})$ is strictly increasing. Assume that $k$ is the minimal natural number such that $k\geq k_{0}$ and $\frac{10^{k^{2}}-k-1}{10^{(k-1)^{2}}}\geq r$. Let $m\geq k$. Since $\vert\sum_{i=n_{m-1}+1}^{n_{m}} \varepsilon_{i} y_{i}\vert\leq 1$ and $\vert y_{i}\vert=\frac{1}{2^{m}}$ for $i\in(n_{m-1},n_{m}]$, we obtain 
$$
\vert\sum\limits_{i\in(n_{m-1},n_{m}]\cap 2\mathbb{N}} \varepsilon_{i}-\sum\limits_{i\in(n_{m-1},n_{m}]\cap (2\mathbb{N}-1)} \varepsilon_{i}\vert\leq 2^{m}.
$$ 
The last inequality means that the excess of ones in the sequence $(\varepsilon_{i})_{i\in(n_{m-1},n_{m}]}$ with odd indexes over those with even indexes is less then $2^{m}$ and vice versa. Consequently, $\vert\sum_{i=n_{m-1}+1}^{n_{m}} \varepsilon_{i} x_{i}\vert\leq \frac{2^{m}}{2^{10^{m^{2}}}}$. Moreover $\sum_{i=1}^{n_{k-1}} \varepsilon_{i} x_{i}=\frac{p_{0}}{2^{10^{(k-1)^{2}}}}$ for some $p_{0}\in\mathbb{Z}$. We have 
$$
\vert\sum_{i=n_{k-1}+1}^{\infty} \varepsilon_{i} x_{i}\vert\leq \sum_{m=k}^{\infty}\vert\sum_{i=n_{m-1}+1}^{n_{m}} \varepsilon_{i} x_{i}\vert\leq\sum_{m=k}^{\infty}\frac{2^{m}}{2^{10^{m^{2}}}}=\sum_{m=k}^{\infty}\frac{1}{2^{10^{m^{2}}-m}}.
$$
Note that
$$
\frac{\frac{1}{2^{10^{(m+1)^{2}}-m-1}}}{\frac{1}{2^{10^{m^{2}}-m}}}=\frac{2^{10^{m^{2}}-m}}{2^{10^{(m+1)^{2}}-m-1}}=2^{10^{m^{2}}-m-10^{(m+1)^{2}}+m+1}=2^{1-(1-10^{2m+1})10^{m^{2}}}\leq\frac{1}{2}
$$
for every $m\in\mathbb{N}$. Hence 
$$
\sum_{m=k}^{\infty}\frac{1}{2^{10^{m^{2}}-m}}\leq\sum_{m=k}^{\infty}\frac{1}{2^{10^{k^{2}}-k}}\cdot\frac{1}{2^{m-k}}=\frac{2}{2^{10^{k^{2}}-k}}=2^{1+k-10^{k^{2}}}.
$$
Since $1+k-10^{k^{2}}\leq -10^{(k-1)^{2}}r$, then
$$\vert\sum\limits_{i=n_{k-1}+1}^{\infty} \varepsilon_{i} x_{i}\vert\leq 2^{-10^{(k-1)^{2}}\cdot r}=\frac{1}{(2^{10^{(k-1)^{2}}})^{r}}.$$
Hence $\vert x-\sum_{i=1}^{n_{k-1}} \varepsilon_{i} x_{i}\vert=\vert\sum_{i=n_{k-1}+1}^{\infty} \varepsilon_{i} x_{i}\vert\leq \frac{1}{(2^{10^{(k-1)^{2}}})^{r}}$. 
Thus $\vert x-\frac{p_{0}}{q_{0}}\vert\leq\frac{1}{q_{0}^{r}}$ with $q_{0}=2^{10^{(k-1)^{2}}}$. That means that either $x$ is a rational number or $x$ is a Liouville number. Finally we obtain $\on{A}(x_{n},y_{n})\subseteq (L\cup \mathbb{Q})\times\mathbb{R}$. Therefore $\on{A}(x_{n},y_{n})$ is of measure zero.
}
\end{example}

\section{Openess of achievement sets}

In this section we show that for some series on the plane its achievement set can be an open set not equal to the whole plane or an open set with two additional points. Such sets are unbounded, since bounded achievement sets are compact. 

\begin{theorem}\label{AlmostOpenAchievementSet}
Let $\sum_{k=1}^{\infty} x_{k}=X<\infty$ with $x_{k}>0$ for every $k\in\mathbb{N}$ and let it satisfy the following conditions: 
\begin{itemize}
\item[(i)] $\on{A}(x_{k})+\on{A}(x_{k})=[0,2X]$
\item[(ii)] for every $a\in(0,2X)$ there exists an interval $I_{a}\subseteq \on{A}(x_{k})$ such that for all $t\in I_{a}$ there exists $z\in \on{A}(x_{k})$ for which $t+z=a$
\end{itemize}
If $(y_{k})$ is conditionally convergent and $\sum_{k=1}^{\infty} y_{k}=Y$ then $\on{A}(\overline{x}_{k},\overline{y}_{k})=(0,2X)\times\mathbb{R}\cup\{(0,0),(2X,Y)\}:=B$, where $\overline{x}_{2k-1}=\overline{x}_{2k}=x_{k}$  and $\overline{y}_{2k-1}=y_{k}, \overline{y}_{2k}=0$ for every $k\in\mathbb{N}$.
\end{theorem}
\begin{proof}
Observe that $\on{A}(\overline{x}_{k},\overline{y}_{k})\subseteq B$. Indeed $\sum_{k=1}^{\infty} (\overline{x}_{k},\overline{y}_{k})=(2\sum_{k=1}^{\infty}x_{k}, \sum_{k=1}^{\infty}y_{k})=(2X,Y)\in \on{A}(\overline{x}_{k},\overline{y}_{k})$ and $(0,0)\in \on{A}(\overline{x}_{k},\overline{y}_{k})$. Moreover $\sum_{k=1}^{\infty}\varepsilon_{k} \overline{x}_{k}\in (0,2X)$ if at the sequence $(\varepsilon_{k})$ there is at least one $0$ and $1$. Hence $\sum_{k=1}^{\infty}\varepsilon_{k} (\overline{x}_{k},\overline{y}_{k})\subseteq B$ for every $(\varepsilon_{k})\in\{0,1\}^{\mathbb{N}}$.

To prove the reversed inclusion, it is sufficient to show that for every $(a,b)\in(0,2X)\times\mathbb{R}$ we have $(a,b)\in  \on{A}(\overline{x}_{k},\overline{y}_{k})$. Let $a\in(0,2X)$ and $b\in\mathbb{R}$. Let $I_{a}\subseteq \on{A}(x_{k})$ be an interval, which satisfies (ii). Then from the absolute convergence of $\sum_{k=1}^{\infty} x_{k}$ one can fix $(\varepsilon_{1}^{a},\ldots, \varepsilon_{k_{a}}^{a})\in\{0,1\}^{k_{a}}$ for which $\sum_{k=1}^{\infty}\varepsilon_{k}^{a}x_{k}\in I_{a}$ for every $(\varepsilon_{k}^{a})_{k=k_{a}+1}^{\infty}\in\{0,1\}^{\mathbb{N}}$. From the conditional convergence of the series  $\sum_{k=1}^{\infty} y_{k}$ we can get the equality $b=\sum_{k=1}^{\infty} \delta_{k} y_{k}$ for some $(\delta_{k})\in\{0,1\}^{\mathbb{N}}$, where $\delta_{k}=\varepsilon_{k}^{a}$ for $k\leq k_{a}$. Then for $t=\sum_{k=1}^{\infty} \delta_{k} x_{k}\in I_{a}$ we define $z=\sum_{k=1}^{\infty} \overline{\delta}_{k} x_{k}$, where $(\overline{\delta}_{k})\in\{0,1\}^{\mathbb{N}}$, such that $t+z=a$. Define the sequence $(\gamma_{k})\in\{0,1\}^{\mathbb{N}}$ by alternating the previous sequences $(\delta_{k})$ and $(\overline{\delta}_{k})$, more precisely $\gamma_{2k-1}=\delta_{k}$ and $\gamma_{2k}=\overline{\delta}_{k}$ for every $k\in\mathbb{N}$. Hence we have: 
$$\sum\limits_{k=1}^{\infty} \gamma_{k} (\overline{x}_{k},\overline{y}_{k})=(\sum\limits_{k=1}^{\infty} \gamma_{k}\overline{x}_{k},\sum\limits_{k=1}^{\infty} \gamma_{k}\overline{y}_{k})=(\sum\limits_{k=1}^{\infty} \delta_{k}x_{k}+\sum\limits_{k=1}^{\infty} \overline{\delta}_{k}x_{k},\sum\limits_{k=1}^{\infty} \delta_{k}y_{k}+\sum\limits_{k=1}^{\infty} \overline{\delta}_{k}\cdot 0)=(t+z,b+0)=(a,b).$$
Finally we get $\on{A}(\overline{x}_{k},\overline{y}_{k})=(0,2X)\times\mathbb{R}\cup\{(0,0),(2X,Y)\}$.
\end{proof}

Now we will construct the series on the plane with an open achievement set. 
\begin{theorem}\label{OpenAchievementSet}
Let a series $\sum_{k=1}^{\infty} v_{k}=X<\infty$ with $v_{k}>0$ for every $k\in\mathbb{N}$ and let it satisfy conditions (i) and (ii) from Theorem \ref{AlmostOpenAchievementSet}. Let $(w_{k})$ be the alternating sequence of non-zero terms of a conditionally convergent series $\sum_{k=1}^{\infty} w_{k}$. If $x_{4k-3}=x_{4k-2}=v_{k}$, $x_{4k-1}=x_{4k}=-v_{k}$ and $y_{2k-1}=w_{k}$, $y_{2k}=0$ for every $k\in\mathbb{N}$ then $\on{A}(x_{k},y_{k})=(-2X,2X)\times\mathbb{R}$.
\end{theorem}

\begin{proof}
The inclusion $\on{A}(x_{k},y_{k})\subset [-2X,2X]\times\mathbb{R}$ is obvious, since $\sum_{k=1}^{\infty} \varepsilon_{k} x_{k}\in [-2 \sum_{k=1}^{\infty} v_{k}, 2\sum_{k=1}^{\infty} v_{k}]$ for every $(\varepsilon_{k})\in\{0,1\}^{\mathbb{N}}$. Moreover, the only way to obtain $2X$ on the first coordinate is to sum up all $x_{k}$'s which are bigger than $0$, more precisely $2X=\sum_{k\in P} x_{k}$ where $P=\{k: x_{k}>0\}=\{4n-i : n\in\mathbb{N} , i\in\{2,3\}\}$.  Since $(w_{k})$ is alternating, we have $\sum_{k\in P} y_{k}= \sum_{k\in\mathbb{N}} w_{2k-1}=\pm\infty$ where the sign $\pm$ depends on the sign of $w_{1}$. Therefore $\on{A}(x_{k},y_{k})\cap\{(2X,y) : y\in\mathbb{R}\}=\emptyset$. In the same way we prove $\on{A}(x_{k},y_{k})\cap\{(-2X,y) : y\in\mathbb{R}\}=\emptyset$ by considering the set of indexes $\mathbb{N}\setminus P$. Hence  $\on{A}(x_{k},y_{k})\subset (-2X,2X)\times\mathbb{R}$.

Now we prove the reverse inclusion. Fix $(a,r)\in (0,2X)\times\mathbb{R}$. Let $I_{a}$ satisfy the assumptions of Theorem \ref{AlmostOpenAchievementSet}. Then we can find an even number $N$ and $(\overline{\varepsilon_{k}})_{k=1}^{N}\in\{0,1\}^N$ such that  $\sum_{k=1}^{N}\overline{\varepsilon_{k}}v_{k}+\sum_{k=N+1}^{\infty}\varepsilon_{k} v_{k}\in I_{a}$ for each $(\varepsilon_{k})_{k=N+1}^{\infty}\in\{0,1\}^\mathbb{N}$. Let $r=\sum_{k=1}^{N}\overline{\varepsilon_{k}}w_{k}+\sum_{k=N+1}^{\infty}\delta_{k} w_{k}$, where $\delta_{k}\in\{0,1\}$ for $k>N$. Define $t=\sum_{k=1}^{N}\overline{\varepsilon_{k}}v_{k}+\sum_{k=N+1}^{\infty}\delta_{k} v_{k}\in I_{a}$.
Then we can find $z=\sum_{k=1}^{\infty} \gamma_{k}v_{k}\in \on{A}(v_{k})$, such that $t+z=a$, where $\gamma_{k}\in\{0,1\}$ for every $k\in\mathbb{N}$. Define $(\alpha_{k})\in\{0,1\}^{\mathbb{N}}$ as follows:
$\alpha_{k}=\overline{\varepsilon_{i}}$ for odd $k\in\{1,\ldots,2N-1\}$, $\alpha_{k}=\delta_{i}$  for odd $k>2N-1$, where $k=2i-1$ and $\alpha_{k}= \gamma_{\frac{k}{2}}$ for even indexes. Hence
$$\sum\limits_{k=1}^{\infty} \alpha_{k} x_{k}=\sum\limits_{k=1}^{N}\overline{\varepsilon_{k}}v_{k}+\sum\limits_{k=N+1}^{\infty}\delta_{k} v_{k}+\sum\limits_{k=1}^{\infty} \gamma_{k}v_{k}=t+z=a$$
and 
$$\sum\limits_{k=1}^{\infty} \alpha_{k} y_{k}=\sum\limits_{k=1}^{N}\overline{\varepsilon_{k}}w_{k}+\sum\limits_{k=N+1}^{\infty}\delta_{k} w_{k}+0=r.$$
Hence  $(0,2X)\times\mathbb{R}\subset \on{A}(x_{k},y_{k})$.
In the similar way we prove that $(-2X,0)\times\mathbb{R}\subset \on{A}(x_{k},y_{k})$. Moreover for a given $r\in\mathbb{R}$ we have $r=\sum_{k=1}^{\infty}\beta_{2k-1} y_{2k-1}$ for some sequence $(\beta_{2k-1})_{k=1}^{\infty}\in\{0,1\}^{\mathbb{N}}$. Define $\beta_{4k}=\beta_{4k-3}$ and $\beta_{4k-2}=\beta_{4k-1}$ for every $k\in\mathbb{N}$. We have 
$$\sum\limits_{k=1}^{\infty}\beta_{k}y_{k}=\sum\limits_{k=1}^{\infty}\beta_{4k-3}y_{4k-3}+
\sum\limits_{k=1}^{\infty}\beta_{4k-2}y_{4k-2}+\sum\limits_{k=1}^{\infty}\beta_{4k-1}y_{4k-1}+
\sum\limits_{k=1}^{\infty}\beta_{4k}y_{4k}=\sum\limits_{k=1}^{\infty}\beta_{2k-1} y_{2k-1}=r$$
and 
$$\sum\limits_{k=1}^{\infty}\beta_{k}x_{k}=\sum\limits_{k=1}^{\infty}\beta_{4k-3}x_{4k-3}+
\sum\limits_{k=1}^{\infty}\beta_{4k-2}x_{4k-2}+\sum\limits_{k=1}^{\infty}\beta_{4k-1}x_{4k-1}+
\sum\limits_{k=1}^{\infty}\beta_{4k}x_{4k}=$$
$$\sum\limits_{k=1}^{\infty}\beta_{4k-3}v_{k}+
\sum\limits_{k=1}^{\infty}\beta_{4k-1}v_{k}+\sum\limits_{k=1}^{\infty}\beta_{4k-1}(-v_{k})+
\sum\limits_{k=1}^{\infty}\beta_{4k-3}(-v_{k})=0$$
We get  $\{0\}\times\mathbb{R}\subset \on{A}(x_{k},y_{k})$. Hence $\on{A}(x_{k},y_{k})=(-2X,2X)\times\mathbb{R}$.
\end{proof}
\begin{example}\emph{
Let $(x_{n},y_{n})_{n\in\mathbb{N}}=((\frac{1}{2},1),(\frac{1}{2},0),(\frac{1}{4},-\frac{1}{2}),(\frac{1}{4},0),(\frac{1}{8},\frac{1}{3}),(\frac{1}{8},0),(\frac{1}{16},-\frac{1}{4}),(\frac{1}{16},0),\ldots)$. Then by Theorem \ref{AlmostOpenAchievementSet} we have $\on{A}(x_{n},y_{n})=(0,2)\times\mathbb{R}\cup\{(0,0),(2,\ln 2)\}$.
}
\end{example}
\begin{example}\emph{
Let $(x_{n},y_{n})_{n\in\mathbb{N}}=((\frac{1}{2},1),(\frac{1}{2},0),(-\frac{1}{2},-\frac{1}{2}),(-\frac{1}{2},0),(\frac{1}{4},\frac{1}{3}),(\frac{1}{4},0),(-\frac{1}{4},-\frac{1}{4}),(-\frac{1}{4},0),\ldots)$. Then by Theorem \ref{OpenAchievementSet} we have $\on{A}(x_{n},y_{n})=(-2,2)\times\mathbb{R}$.
}
\end{example}

\section{Achievement sets of potentially conditionally convergent series}

In this section we consider the set $\bigcup_{\sigma\in S_{\infty}}\on{A}(x_{\sigma(k)})$ which is the union of achievement sets of all rearrangements of a series $\sum_{k=1}^\infty x_k$. If the underlay series is absolutely convergent, then $\bigcup_{\sigma\in S_{\infty}}\on{A}(x_{\sigma(k)})$ equals $\on{A}(x_n)$. Here we study a situation when the given series is conditionally convergent. Clearly only the multidimensional case is interesting. 

We say that the series $\sum_{k=1}^{\infty} x_{k}$ is potentially conditionally convergent if there exists a permutation $\sigma$ for which the rearranged series  $\sum_{k=1}^{\infty} x_{\sigma(k)}$ is conditionally convergent. For example on the real line $\sum_{k=1}^{\infty} x_{k}$ is potentially conditionally convergent if and only if $x_{n}\rightarrow 0$ and $\sum_{\{k :x_{k}>0\}} x_{k}=\infty$, $\sum_{\{k :x_{k}<0\}} x_{k}=-\infty$. We will also consider the set $\on{A}_{\on{abs}}(x_{k})=\{\sum_{k=1}^{\infty} \varepsilon_{k} x_{k}  : \sum_{k=1}^{\infty} \varepsilon_{k}\Vert x_{k}  \Vert<\infty,  \varepsilon_{k}\in\{0,1\} \ \text{for each} \  k\in\mathbb{N} \}$.

\begin{theorem}\label{TheoremRearrangedAchievementSet}
Let $\sum_{k=1}^{\infty} y_{k}$ be an absolutely convergent series such that the function $(\varepsilon_{k})\rightarrow \sum_{k=1}^{\infty} \varepsilon_{k} y_{k}$, where $(\varepsilon_{k})\in \{0,1\}^{\mathbb{N}}$ is injective. Assume that $\sum_{k=1}^{\infty} x_{k}$ is conditionally convergent. Let $X=\{\sum_{k=1}^\infty\varepsilon_kx_k:\sum_{k=1}^\infty\varepsilon_kx_k$ is potentially conditionally convergent$\}$. Then
\begin{itemize}
\item[(i)] $\bigcap_{\sigma\in S_{\infty}} \on{A}( x_{\sigma(k)}, y_{\sigma(k)})=\on{A}_{\on{abs}}(x_{k},y_{k})$.
\item[(ii)]  $\bigcup_{\sigma\in S_{\infty}} \on{A}( x_{\sigma(k)}, y_{\sigma(k)})=(X\times\R)\cup \on{A}_{\on{abs}}(x_{k},y_{k})$.
\end{itemize}
\begin{proof}
Let $(x,y)\in \on{A}_{\on{abs}}(x_{k},y_{k})$ and $\sigma\in S_{\infty}$, then $(x,y)=\sum_{k=1}^{\infty}\varepsilon_{k}(x_{k}, y_{k})=\sum_{k=1}^{\infty}\varepsilon_{\sigma(k)}(x_{\sigma(k)}, y_{\sigma(k)})$, because $\sum_{k=1}^{\infty}\varepsilon_{k}\Vert(x_{k}, y_{k})\Vert=\sum_{k=1}^{\infty}\Vert\varepsilon_{k}(x_{k}, y_{k})\Vert<\infty$ and we know that every absolutely convergent series is unconditionally convergent. Hence $(x,y)\in \bigcap_{\sigma\in S_{\infty}} \on{A}( x_{\sigma(k)}, y_{\sigma(k)})$.

Suppose that $(x,y)\notin \on{A}_{\on{abs}}(x_{k},y_{k})$. If $(x,y)\notin \on{A}(x_{k},y_{k})$, then $(x,y)\notin \bigcap_{\sigma\in S_{\infty}} \on{A}( x_{\sigma(k)}, y_{\sigma(k)})$. Assume that $(x,y)\in \on{A}(x_{k},y_{k})$ and $(x,y)=\sum_{k=1}^{\infty}\varepsilon_{k}(x_{k}, y_{k})$. Since $(x,y)\notin \on{A}_{\on{abs}}(x_{k},y_{k})$, we have $\sum_{k=1}^{\infty}\varepsilon_{k}\Vert(x_{k}, y_{k})\Vert=\infty$. One can find $\sigma\in S_{\infty}$ such that $\sum_{k=1}^{\infty}\varepsilon_{\sigma(k)}x_{\sigma(k)}=\overline{x}\neq x$. From the absolute convergence of $\sum_{k=1}^{\infty} y_{k}$, we have $\sum_{k=1}^{\infty}\varepsilon_{\sigma(k)}(x_{\sigma(k)},y_{\sigma(k)})=(\overline{x},y)$. We assumed that the second coordinate is obtained for at most one sequence $(\varepsilon_{k})$, so $\sigma = \on{id}$ and in consequence $\overline{x}=x$, which give us contradiction.

Now we will prove the second equality. 
\\,,$\subseteq$''. Let $(x,y)=\sum_{k=1}^{\infty}\varepsilon_{k}(x_{\sigma(k)}, y_{\sigma(k)})$ for some $(\varepsilon_{k})\in\{0,1\}^{\mathbb{N}}$ and $\sigma\in S_{\infty}$. We have two posibilities: 
\begin{enumerate}
\item If the series $\sum_{k=1}^{\infty}\varepsilon_{k}y_{\sigma(k)}$ is absolutely convergent, then $\sum_{k=1}^{\infty} \varepsilon_{\sigma^{-1}(k)}x_{k}=\sum_{k=1}^{\infty}\varepsilon_{k}x_{\sigma(k)}=x$ and also $\sum_{k=1}^{\infty} \varepsilon_{\sigma^{-1}(k)}y_{k}=\sum_{k=1}^{\infty}\varepsilon_{k}y_{\sigma(k)}=y$. Hence the series $\sum_{k=1}^{\infty} \varepsilon_{\sigma^{-1}(k)}(x_k,y_{k})$ is absolutely convergent, so $(x,y)\in \on{A}_{\on{abs}}(x_{k},y_{k})$.
\item If the series $\sum_{k=1}^{\infty}\varepsilon_{k}y_{\sigma(k)}$ is conditionally convergent, then it is also potentially conditionally convergent, so $(x,y)\in X\times\mathbb{R}$.
\end{enumerate}
,,$\supseteq$''. From the first equality we know that $\on{A}_{\on{abs}}(x_{k},y_{k})=\bigcap_{\sigma\in S_{\infty}} \on{A}( x_{\sigma(k)}, y_{\sigma(k)})\subset \bigcup_{\sigma\in S_{\infty}} \on{A}( x_{\sigma(k)}, y_{\sigma(k)})$, so it is enough to show that $\mathbb{R}\times Y\subset \bigcup_{\sigma\in S_{\infty}} \on{A}( x_{\sigma(k)}, y_{\sigma(k)})$. Fix $a\in\mathbb{R}$ and $y\in Y$, then there exist a sequence $(\varepsilon_{k})\in\{0,1\}^{\mathbb{N}}$ and $\sigma\in S_{\infty}$ such that $y=\sum_{k=1}^{\infty}\varepsilon_{k} y_{\sigma(k)}$ and $\sum_{k=1}^{\infty}\varepsilon_{k} x_{\sigma(k)}$ converge. One can rearrange the terms of the conditionally convergent series $\sum_{n=1}^\infty x_n$ such that  $a=\sum_{k=1}^{\infty}\varepsilon_{\tau(k)} x_{\tau(\sigma(k))}$. Since $\sum_{k=1}^{\infty} y_{k}$ is absolutely convergent, then $\sum_{k=1}^{\infty}\varepsilon_{\tau(k)} y_{\tau(\sigma(k))}=\sum_{k=1}^{\infty}\varepsilon_{k} y_{\sigma(k)}=y$. Hence $(a,y)\in \on{A}(x_{\tau(\sigma(k))},y_{\tau(\sigma(k))})$.
\end{proof}
\end{theorem}

By Theorem \ref{AIsDenseInProduct} the achievement set $\on{A}(x_n,y_n)$ is dense in $\R\times\on{A}(y_n)$. So is $S=\bigcup_{\sigma\in S_{\infty}} \on{A}( x_{\sigma(k)}, y_{\sigma(k)})$. By Theorem \ref{TheoremRearrangedAchievementSet} a vertical section $S_x$ equals to $\R$ if $x\in X$; it is a singleton if $x=\sum_{k=1}^\infty\varepsilon_kx_k$ is absolutely convergent, and it is empty if $\sum_{k=1}^\infty\varepsilon_kx_k$ is not potentially conditionally convergent. Hence the series $\sum_{k=1}^\infty\varepsilon_kx_k$ has only finitely many positive elements or only finitely many negative elements.  
\vspace{0.5 cm}

Let us finish the paper with the list of open questions:
\begin{enumerate}
\item Does there exist a conditionally convergent series $\sum_{n=1}^{\infty} x_n$ on the plane such that $\on{A}(x_n)$ is a graph of a function with a domain being a bounded interval?
\item Let $\on{SR}(x_n)=\R^k$. Is it true that either $\on{A}(x_n)=\R^k$ or $\on{A}(x_n)$ is of measure zero? 
\item Lemma \ref{AIsClosedInSR} implies that the achievement sets of conditionally convergent series in finite dimensional spaces are unbounded. On the other hand there is an example of conditionally convergent series in $c_0$ with a closed and bounded achievement set, see Example \ref{C0Example}. Obviously such series can be found in every Banach space containing isomorphic copy of $c_0$. Note that the series from Example \ref{C0Example} has the unbounded achievement set in $\ell_1$ (and it is a well-known fact that $\ell_1$ does not contain a copy of $c_0$). Is there a conditionally convergent series in $\ell_1$ with a bounded and closed achievement set?
\item In the proof of Theorem \ref{NeitherFSigmaNorGDelta}, we show that the achievement set is analytic. Moreover there are analytic sets which are not Borel. Is there a (conditionally convergent) series which achievement set is non-Borel?  
\end{enumerate}

\end{document}